\documentclass[11pt,final]{amsart}
\usepackage[utf8]{inputenc}
\usepackage[T1]{fontenc}
\usepackage[letterpaper,centering]{geometry}
\usepackage{lmodern}
\usepackage{eucal}
\usepackage{mathrsfs}
\usepackage{amsmath,amssymb,amsfonts,amscd}
\usepackage{microtype,xcolor,enumitem,booktabs,needspace}
\definecolor{darkblue}{rgb}{0,0,0.3}
\definecolor{darkgreen}{rgb}{0,0.4,0}
\usepackage[colorlinks=true,linkcolor=darkblue,citecolor=darkgreen,
 urlcolor=darkblue,unicode,pdfborder={0 0 0}]{hyperref}
\hypersetup{pdftitle={Campana strong approximation on quadrics with hyperplane boundary},
 pdfsubject={Campana points, quadrics and strong approximation}}
\setlist[enumerate]{label={\upshape(\arabic*)},topsep=.7ex,leftmargin=*}
\setlist[itemize]{leftmargin=*}
\theoremstyle{plain}
\newtheorem{thm}{Theorem}[section]
\newtheorem{lem}[thm]{Lemma}
\newtheorem{prop}[thm]{Proposition}
\newtheorem{cor}[thm]{Corollary}
\theoremstyle{definition}
\newtheorem{defn}[thm]{Definition}
\newtheorem{rmk}[thm]{Remark}
\numberwithin{equation}{section}
\newcommand{\PP}{\mathbf P}
\newcommand{\Aff}{\mathbf A}
\newcommand{\A}{\mathbf A}
\newcommand{\OO}{\mathcal O}
\newcommand{\FF}{\mathbf F}
\newcommand{\ZZ}{\mathbf Z}
\newcommand{\kbar}{\overline{k}}
\newcommand{\sX}{\mathscr X}
\newcommand{\sY}{\mathscr Y}

\newcommand{\sD}{\mathscr D}

\newcommand{\sQ}{\mathscr Q}

\newcommand{\CC}{\mathcal C}

\newcommand{\sm}{\mathrm{sm}}
\DeclareMathOperator{\Spec}{Spec}
\DeclareMathOperator{\Supp}{Supp}
\DeclareMathOperator{\ord}{ord}

\newcommand{\N}[2]{N_{#1/#2}}
\title[Campana points on quadrics]
{Campana points on quadrics with hyperplane boundary}
\author{Dasheng Wei}
\address{Academy of Mathematics and Systems Science, CAS, Beijing 100190,
	P.\ R.\ China {\normalfont\itshape and} School of Mathematical Sciences,
	University of CAS, Beijing 100049, P.\ R.\ China}
\email{dshwei@amss.ac.cn}
\date{}

\begin{document}

\begin{abstract}
Let $Q\subset\PP^{d+1}_k$ be a smooth quadric of dimension $d\ge2$
over a number field, and let $H_i=Q\cap\Pi_i$ be hyperplane
sections for $1\le i\le d$.
For every integer $m\ge2$, we prove Campana strong approximation
for the pair
$\bigl(Q,(1-1/m)(H_1+\cdots+H_d)_{\mathrm{red}}\bigr)$.
\end{abstract}

\maketitle

\section{Introduction}


Quadrics form a classical family in number theory.
By the Hasse--Minkowski theorem, they satisfy the Hasse principle
and weak approximation.
For affine quadrics of dimension at least $3$, strong approximation
holds under a non-compactness condition;
in dimension $2$, it holds with the Brauer--Manin obstruction
\cite[Theorem~8.3]{CTX}.
The Hasse principle and weak approximation for intersections of two
quadrics have also been extensively investigated; see, for instance,
\cite{CTSSD,CTSSDII,SS,Wittenberg,HeathBrown,CT}.

In the Campana setting, Nakahara and Streeter established
approximation results for projective spaces, including the projective
line with a boundary of degree two \cite{NS}. The log Fano case on
$\PP^1$ is treated in \cite[Theorem~1.1]{WeiProjective}.
Mitankin, Nakahara, and
Streeter proved Campana strong approximation for a smooth quadric
whose boundary is a single hyperplane section
\cite[Theorem~1.4]{MNS}. In this paper, we treat the case where
the boundary is the reduced union of hyperplane sections.

For $d\ge2$, let $Q\subset\PP^{d+1}_k$ be a smooth quadric, let
$\Pi_1,\ldots,\Pi_d$ be hyperplanes defined over $k$, and put
\begin{equation}\label{eq:pair}
	H_i=Q\cap\Pi_i,\qquad D=H_1+\cdots+H_d,
	\qquad \Delta_m=\left(1-\frac{1}m\right)D_{\mathrm{red}}.
\end{equation}
Each distinct prime component of $D$ occurs in
$D_{\mathrm{red}}$ with multiplicity one. The hyperplanes need
not be distinct. 
For $d\ge2$, every hyperplane section is geometrically reduced,
since its quadratic form has rank at least $d$. For $d\ge3$,
every hyperplane section is geometrically integral. On a surface,
distinct sections may share a line.

\begin{thm}\label{thm:main}
Let $Q\subset\PP^{d+1}_k$ be a smooth quadric of dimension
$d\ge2$. Let
$\Pi_1,\ldots,\Pi_d$ be hyperplanes defined over $k$, and let
$D$ and $\Delta_m$ be as in \eqref{eq:pair}.
Then, for every integer $m\ge2$, the pair $(Q,\Delta_m)$
satisfies Campana strong approximation, in the sense that
\[
 \overline{(Q\setminus D)(k)}
 =(Q,\Delta_m)(\A_k).
\]
Here the closure is taken in the Campana topology
of Definition~\ref{def:strong}.
\end{thm}

As an immediate consequence of Theorem~\ref{thm:main}, we obtain the following corollary concerning $S$-Campana points.
\begin{cor}\label{cor:original-notation}
Under the hypotheses of Theorem~\ref{thm:main}, fix a finite
$S\supseteq\Omega_k^\infty$ and a proper flat integral
$\OO_S$-model $\sQ$ of $(Q,\Delta_m)$. The strict $S$-Campana points
of this model are dense in
\[
 \mathcal A_S(\sQ,\Delta_m)
 =\prod_{v\in S}Q(k_v)\times
 \prod_{v\notin S}\CC_v(\sQ,\Delta_m)
\]
with its product topology, where $\CC_v$ is defined using this
fixed model. 
\end{cor}

By monotonicity of Campana strong approximation \cite[Lemma 2.5(1)]{WeiMarkoff}, we immediately obtain the following result for unequal weights.
\begin{cor}\label{cor:unequal-weights}
	Under the hypotheses of Theorem~\ref{thm:main}, write
	$D_{\mathrm{red}}=\sum_{\alpha\in I}D_\alpha$ as the sum of its distinct
	prime components. Let $n_\alpha\ge2$ be arbitrary integers and put
	\[
		\Delta_{\boldsymbol n}
		=\sum_{\alpha\in I}\left(1-\frac1{n_\alpha}\right)D_\alpha.
	\]
	Then $(Q,\Delta_{\boldsymbol n})$ satisfies Campana strong
	approximation. 
\end{cor}

%

The common intersection
$$E=Q\cap\Pi_1\cap\cdots\cap\Pi_d$$ is understood
scheme-theoretically.
First suppose that $E(k)\ne\varnothing$, allowing $E$ to have
positive dimension. Let $L=\Pi_1\cap\cdots\cap\Pi_d$.
Choose a line $\ell\subset L$ through a rational point of $E$.
Let $\Pi$ be a general plane containing $\ell$. If $\ell\subset Q$, then
$\Pi\cap Q=\ell\cup C$,
where $C$ is the residual line, whose reduced boundary is the single point $x=C\cap\ell$. Otherwise,
$C=\Pi\cap Q$ is a smooth conic with reduced boundary $(C\cap\ell)_{\mathrm{red}}$ of degree at most two.
Proposition~\ref{prop:rational-base-locus} treats all these cases.

We may therefore assume $E(k)=\varnothing$.
When $E$ is finite, the linear forms defining the hyperplanes
$\Pi_i$ are linearly independent. Then
$\Lambda=\Pi_1\cap\cdots\cap\Pi_d\simeq\PP^1_k$ is not
contained in $Q$, $D$ is reduced, and $E$ has length two.
Since $E(k)=\varnothing$, it is a closed point of degree two. 
For surfaces, we use the conic pencil through
$H_1\cap H_2$. 
Using the fibration method of \cite{HWW}, we choose a smooth fibre
with local points in the prescribed Campana-open neighborhood.
Campana strong approximation on $\PP^1$ then gives the
desired rational point. For $d>2$, the common plane conic
is automatically reduced.
We may choose a suitable quadric surface through
$Q\cap\Pi_2\cap\cdots\cap\Pi_d$ and thereby reduce the
approximation problem to the surface case. 

It remains to treat positive-dimensional $E$ with
$E(k)=\varnothing$. 
Choose a plane $L'\subset L$ with smooth
conic $B=Q\cap L'$. For general $R$, the section
$\Sigma_R=Q\cap\langle L',R\rangle$ is a smooth quadric
surface, and every $H_i$ restricts to $B$.
The one-boundary surface case in
Corollary~\ref{cor:surface-one-boundary} completes this reduction.

The paper is organised as follows. After the preliminaries in
Section~\ref{sec:prelim},
Section~\ref{sec:rational-base-locus} treats all $E$ with
$E(k)\ne\varnothing$. The remaining surface case, which is the
hardest part of the paper, is proved in Section~\ref{sec:surface}.
Section~\ref{sec:higher} carries out the remaining reductions to
surfaces in both the finite and positive-dimensional cases.

\noindent\textbf{Acknowledgements.}
This work is supported by the National Natural Science Foundation
of China (Grant Nos.~12371014 and 12231009).

\section{Campana approximation and preliminary results}\label{sec:prelim}

Let $\Omega_k$ be the set of places of $k$, and let
$\Omega_k^\infty$ be the subset of archimedean places. For a finite place~$v$, we write $\OO_v$, $\mathfrak m_v$ and $\FF_v$ for its valuation
ring, maximal ideal and residue field. We normalize $\ord_v$ to
have image $\ZZ$.
Throughout this paper, we assume that $S\supseteq\Omega_k^\infty$ is a finite set of places. Let $\OO_S$ be the ring of
$S$-integers.

\subsection{Campana topology}\label{subsec:local-conditions}

Let $X$ be a smooth geometrically integral $k$-variety, not
necessarily proper. Let $D=\sum_iD_i$ be a reduced effective divisor
with distinct prime components defined over $k$, and put
\[
 \Delta=\sum_i\left(1-\frac1{m_i}\right)D_i,
 \qquad m_i\in\ZZ_{\ge2}.
\]
All weights on $X$ are finite. For equal weights $m_i=m$, this gives $\Delta=(1-1/m)D$. Fix a finite set
$S$ and a flat separated integral
$\OO_S$-model $\sX$ of $X$.
Let $\sD_i$ be the reduced closure of $D_i$.

For $v\notin S$, write
$D_i\times_k k_v=\bigcup_\lambda D_{i,\lambda}$
for the local irreducible components, and let $\sD_{i,\lambda}$
be their reduced closures in $\sX_{\OO_v}$.
For $x\in\sX(\OO_v)\subset X(k_v)$, let
$\widetilde x:\Spec\OO_v\to\sX_{\OO_v}$ be its section.
\needspace{7\baselineskip}
For any closed subscheme $\mathscr Z\subset\sX_{\OO_v}$,
let $\mathcal I_{\mathscr Z}\subset\OO_{\sX_{\OO_v}}$
be its ideal sheaf. Define its pulled-back ideal along
$\widetilde x$ by
\[
 J_{\mathscr Z,x}
 =\widetilde x^{-1}\mathcal I_{\mathscr Z}\cdot\OO_v
 \subseteq\OO_v.
\]
If $f_1,\ldots,f_r$ generate $\mathcal I_{\mathscr Z}$ near
the reduction of $x$, then
$J_{\mathscr Z,x}=(f_1(x),\ldots,f_r(x))$.
In particular, write
$J_{i,\lambda,x}:=J_{\sD_{i,\lambda},x}$.
Since $\OO_v$ is a discrete valuation ring, a nonzero
$J_{i,\lambda,x}$ equals $\mathfrak m_v^{\,n}$ for a unique
integer $n\ge0$. We define the \emph{intersection multiplicity} by
\[
 n_v(\sD_{i,\lambda},x)=
 \begin{cases}
 \operatorname{length}_{\OO_v}(\OO_v/J_{i,\lambda,x}),
       &J_{i,\lambda,x}\ne0,\\
 \infty,&J_{i,\lambda,x}=0.
 \end{cases}
\]

Define
\begin{align*}
 \CC_v(\sX,\Delta)
 &=\{x\in\sX(\OO_v): n_v(\sD_{i,\lambda},x)
             \in\{0\}\cup\ZZ_{\ge m_i}\cup\{\infty\}
             \text{ for every }i,\lambda\},&&v\notin S,\\
 \CC_v(\sX,\Delta)&=X(k_v),&&v\in S.
\end{align*}
We use the local-component convention of
\cite[Definitions~2.2--2.4]{NS}. A local point $x\in\CC_v(\sX,\Delta)$ is called a
\emph{strict local Campana point} if $x\in(X\setminus D)(k_v)$.
A rational point $x\in(X\setminus D)(k)$ is called a
\emph{strict $S$-Campana point} for the fixed model if it
extends to an $\OO_v$-point of $\sX$ for every $v\notin S$
and satisfies
\[
 n_v(\sD_{i,\lambda},x)=0
 \quad\text{or}\quad
 n_v(\sD_{i,\lambda},x)\ge m_i
 \qquad\text{for every }v\notin S,\ i \text{ and }\lambda.
\]

For the fixed $\OO_S$-model $\sX$ and a finite set $T\supseteq S$, write
\[
 \mathcal A_T(\sX,\Delta)=\prod_{v\in T}X(k_v)\times \prod_{v\notin T}\CC_v(\sX,\Delta)
\]
with the product topology.

\begin{defn}[Campana topology]\label{def:strong}
The \emph{Campana adelic space} is
\[
 (X,\Delta)(\A_k)
 =\prod_v'\bigl(X(k_v),\CC_v(\sX,\Delta)\bigr) =\bigcup_{T\supseteq S\ {\rm finite}}\mathcal A_T(\sX,\Delta).
\]
Its points satisfy the local Campana conditions at almost
every place. We equip it with the restricted-product topology,
which we call the \emph{Campana topology}. It is independent of the choice of the model
since any two models
agree outside
finitely many places. A basis consists of
\[
 \prod_{v\in T}W_v\times\prod_{v\notin T}\CC_v(\sX,\Delta),
 \qquad T\supseteq S\text{ finite},
\]
where $W_v\subset X(k_v)$ is open in the usual local
analytic topology.  
Following \cite[Definition~2.1]{WeiMarkoff}, the pair
$(X,\Delta)$ satisfies \emph{Campana strong approximation}
if $(X\setminus D)(k)$ is dense in $(X,\Delta)(\A_k)$.

\end{defn}

For a fixed model $\sX$, \emph{Campana weak approximation}
means that its strict $S$-Campana points are dense in
$\mathcal A_S(\sX,\Delta)$ with the product topology; see \cite[Definition~2.7(ii)]{NS}.


We recall \cite[Lemma~2.3 and Proposition~2.4(1)]{WeiMarkoff}.
\begin{lem}\label{lem:local-open}
	\begin{enumerate}
		\item For finite weights, $\CC_v(\sX,\Delta)$ is open and closed in
		$X(k_v)$. 
		\item Campana strong approximation for $(X,\Delta)$ is
		equivalent to Campana weak
		approximation for every finite $S$ and
		every flat separated integral model.
	\end{enumerate}


\end{lem}

\subsection{Orbifold morphisms and continuity}\label{subsec:orbifold-morphisms}

\begin{defn}[Orbifold morphism]\label{def:orbifold-morphism}
Let $(Y,\Delta_Y)$ and $(X,\Delta_X)$ be smooth
geometrically integral $k$-varieties with finite-weight
Campana boundaries. Each geometric boundary component
inherits its weight.  Let $m_Y(F)=1$ for a geometric
prime divisor $F$ outside the boundary of $Y_{\kbar}$.
Following \cite[Definition~2.3]{Campana}, a $k$-morphism
$f:Y\to X$ is a \emph{morphism of Campana orbifolds} if
$f(Y)\not\subseteq\Supp\Delta_X$ and, for every geometric
boundary component $D$ of $X_{\kbar}$,
\[
 f_{\kbar}^*D=\sum_F e_{F,D}F\quad(e_{F,D}>0)
 \quad\Longrightarrow\quad
 e_{F,D}\,m_Y(F)\ge m_X(D).
\]
Here $f_{\kbar}^*D$ is the effective Cartier pullback of
$D$. The sum runs over the prime components
of its support, and $e_{F,D}$ is the order of vanishing along $F$.
For an inclusion $j:Y\hookrightarrow X$, this
is the scheme-theoretic divisor restriction $D|_Y$, with
its multiplicities.
\end{defn}


\begin{prop}[Continuity]\label{prop:orbifold-continuity}
An orbifold morphism $f:(Y,\Delta_Y)\to(X,\Delta_X)$
induces a continuous map
\[
 f_{\mathbf A}:(Y,\Delta_Y)(\A_k)
          \longrightarrow(X,\Delta_X)(\A_k).
\]
More precisely, for any fixed models $\sY$ of $Y$ and
$\sX$ of $X$, there is a finite set $T\supseteq\Omega_k^\infty$ such that
\begin{equation}\label{inclusion}
 f\bigl(\CC_v(\sY,\Delta_Y)\bigr)
       \subseteq\CC_v(\sX,\Delta_X)\qquad(v\notin T).
\end{equation}
In particular, the inverse image of a Campana-open set is open.
\end{prop}

\begin{proof}
Choose a finite Galois extension $L/k$ over which all
geometric components of the boundaries and of the divisors
$f_{\kbar}^*D$ are defined. 
Choose a finite set $T$ containing
the ramified places
of $L/k$ such that $f$ extends to the models and the pullback identities
hold for the flat Cartier closures over
the integral closure of $\OO_T$ in $L$.
For $v\notin T$ and $w\mid v$, the $k_v$-irreducible
components correspond to $\operatorname{Gal}(L_w/k_v)$-orbits
of the geometric components. 
For a local irreducible component $D_\lambda$ of the boundary on $X$, write
$f^*D_\lambda=\sum_\alpha e_{\alpha,\lambda}F_\alpha$
$(e_{\alpha,\lambda}>0)$.
The geometric inequalities imply
$e_{\alpha,\lambda}m_Y(F_\alpha)\ge m_X(D_\lambda)$.
For
$x\in\CC_v(\sY,\Delta_Y)$ and $v\notin T$, we have
\[
 n_v(\mathscr D_\lambda,f(x))
 =\sum_\alpha
    e_{\alpha,\lambda}\,n_v(\mathscr F_\alpha,x).
\]
Every positive $n_v(\mathscr F_\alpha,x)$ is at least
$m_Y(F_\alpha)$.
Hence every positive $n_v(\mathscr D_\lambda,f(x))$ is
at least $m_X(D_\lambda)$. This proves the inclusion (\ref{inclusion}).
The inclusion shows that $f_{\mathbf A}$ is well-defined.
To prove continuity, let $\mathcal N=\prod_v N_v$ be a
basic neighborhood of $f_{\mathbf A}((y_v))$
in $(X,\Delta_X)(\A_k)$.
Choose a finite $B$ containing $T$, the prescribed places
of $\mathcal N$, and all places where
$y_v\notin\CC_v(\sY,\Delta_Y)$.
At $v\in B$, local continuity gives an open neighborhood
$U_v$ of $y_v$ with $f(U_v)\subseteq N_v$.
Then
\[
 \prod_{v\in B}U_v\times
 \prod_{v\notin B}\CC_v(\sY,\Delta_Y)
 \subseteq f_{\mathbf A}^{-1}(\mathcal N)
\]
is a basic neighborhood of $(y_v)$ in
$(Y,\Delta_Y)(\A_k)$. This proves continuity.
\end{proof}



\subsection{Approximation results}

We use the following result on curves.

\begin{prop}[Nakahara--Streeter]\label{prop:curve}
Let $C\simeq\PP^1_k$, and let $E\subset C$ be a reduced effective
divisor of degree at most two. Assign to each closed point
$P\in\Supp(E)$ an arbitrary integer $n_P\ge2$ and put
\[
 \Delta_E=\sum_{P\in\Supp(E)}\left(1-\frac1{n_P}\right)P.
\]
Then $(C,\Delta_E)$ satisfies Campana strong approximation.
\end{prop}

\begin{proof} 
	It  suffices to show that Campana weak 
	approximation holds for  a given $\OO_S$-model and every finite $T\supset S$.
If $E=0$, use ordinary weak approximation on $\PP^1$.
If every point of $E$ is rational, apply
\cite[Proposition~3.5]{NS} with $n=1$. Otherwise, $E$ is a
single closed point of degree two, and
\cite[Proposition~3.15]{NS} applies. 
Lemma~\ref{lem:local-open}\,(2) gives Campana strong approximation. 
For a more general statement, see
\cite[Theorem~1.1]{WeiProjective}. 
\end{proof}

We also use ordinary strong approximation on affine cones
in the following case.

\begin{thm}[Wei--Xu--Zhu]\label{thm:cone}
Let $V\subset\PP^n_k$, $n\ge7$, be a pure geometrically integral
intersection of two quadrics, not a cone, and suppose that $V$
has finitely many singular points. Let $W$ be its affine cone
and let $Z\subset W$ be closed of codimension at least two.
Then $W^{\sm}\setminus Z$ satisfies strong approximation off
any one place of $k$.
\end{thm}

This is \cite[Theorem~1.6(iii)]{WXZ}. In particular, no
Brauer--Manin condition is required in this dimension.


We follow the local-component
convention of Section~\ref{subsec:local-conditions} and the
Campana topology of Definition~\ref{def:strong}. We use \cite[Theorem~3.17]{MNS} in the following form.

\begin{lem}\label{lem:campana-formal}
Let $Y$ be smooth and geometrically integral, and let $\Delta$
have finite weights.
Let $V\subset Y\setminus\Supp\Delta$ be a dense open subset,
and let $\mathcal B\subset\operatorname{Br}V$ be a finite subgroup.
Suppose $(p_v)\in(Y,\Delta)(\A_k)$ is orthogonal
to $\mathcal B\cap\operatorname{Br}Y$.
For every neighborhood $\mathcal N$ of $(p_v)$ in the Campana
topology, there exists
\[
 (p_v')\in
 (Y,\Delta)(\A_k)\cap V(\A_k)^{\mathcal B}
 \quad\text{with}\quad (p_v')\in\mathcal N.
\]
Thus $(p_v')$ can be chosen arbitrarily close to $(p_v)$
in the Campana topology.
\end{lem}

The nonproper case follows by assigning infinite weights at infinity
in a smooth proper compactification with divisorial boundary.
The usual argument for Harari's formal lemma \cite[Corollary~2.6.1]{Harari94} handles arbitrary $V$,
using the local corrections of \cite[Theorem~3.16]{MNS} to preserve
the componentwise Campana conditions.

\section{The case $E(k)\neq \varnothing$}\label{sec:rational-base-locus}

Let $Q\subset\PP(V)$ and $H_1,\ldots,H_d$ be as in
Theorem~\ref{thm:main}, and let $q$ be the nondegenerate
quadratic form defining $Q$. Write
\[
 L=\bigcap_{i=1}^d\Pi_i=\PP(W),\qquad E=Q\cap L.
\]
We do not assume that $E$ is finite.

Let $\ell\subset L$ be a $k$-line. If $\ell\subset Q$,
then for general $R\in Q\setminus D$ the plane
$\langle\ell,R\rangle$ is not contained in $Q$, and
\[
 Q\cap\langle\ell,R\rangle=\ell\cup C_R,
\]
where $C_R$ is a distinct $k$-line through $R$.
Every $\Pi_i$ cuts this plane along $\ell$, so
$H_i|_{C_R}=P_R$, where $P_R=\ell\cap C_R\in E(k)$.

If $\ell\not\subset Q$, put $Z=Q\cap\ell$.
For general $R\in Q\setminus D$,
the plane $\langle\ell,R\rangle$ cuts out a smooth conic
$C_R$, and $H_i|_{C_R}=Z$ for every $i$.
Indeed, the restriction of $q$ to the two-dimensional space
defining $\ell$ has rank two or one. In rank two, adjoining a
vector of nonzero norm in its orthogonal complement gives a
nondegenerate three-dimensional restriction. The resulting
smooth conic contains a geometric point $R\notin\ell$.
Thus the smoothness locus in $Q\setminus\ell$ is nonempty;
it meets $Q\setminus D$ by irreducibility of $Q$.
In rank one, $Z=2P$ for some $P\in Q(k)$,
and the conic is smooth whenever $R\notin T_PQ$.

\begin{lem}[Reduction to a rational curve]\label{lem:conic-reduction}
Let $C\subset Q$ be a smooth plane conic or a $k$-line,
not contained in $D$, and suppose
$F=(D|_C)_{\mathrm{red}}$ has degree at most two.
Put $\delta_C=(1-1/m)F$ and let $j:C\hookrightarrow Q$.
If a Campana-open set $\mathcal N\subset(Q,\Delta_m)(\A_k)$
has a nonempty inverse image
\[
 \mathcal G=j_{\mathbf A}^{-1}(\mathcal N)
       \subset(C,\delta_C)(\A_k),
\]
then $\mathcal N$ contains a point of $(C\setminus F)(k)$.
\end{lem}

\begin{proof}
Each geometric boundary pullback is supported on $F$, so $j$
is an orbifold morphism and $\mathcal G$ is open by
Proposition~\ref{prop:orbifold-continuity}.
Its nonemptiness gives local points on $C$ everywhere;
the Hasse principle yields $C\simeq\PP^1_k$ when $C$ is a conic.
Proposition~\ref{prop:curve} gives a strict rational point
in $\mathcal G$, whose image lies in $\mathcal N$.
\end{proof}

\begin{prop}\label{prop:rational-base-locus}
Under the hypotheses of Theorem~\ref{thm:main}, assume that
$E(k)\ne\varnothing$. Then $(Q,\Delta_m)$ satisfies
Campana strong approximation for every $m\ge2$.
\end{prop}

\begin{proof}

Choose $A\in E(k)$ and a $k$-line $\ell\subset L$ through $A$. For a given proper model $\sQ$ of $Q$,
let $\mathcal N\subset(Q,\Delta_m)(\A_k)$ be nonempty and open. 
We may assume  $$\mathcal N= \prod_{v\in T}W_v\times\prod_{v\notin T}\CC_v(\sQ,\Delta_m).$$

Suppose first that $\ell\not\subset Q$. We enlarge $T$ so that $A$ satisfies the local Campana condition outside $T$. Set
$W_v=\CC_v(\mathscr Q,\Delta_m)$
at every newly added place. Using weak approximation on $Q$, we may choose a general 
$R\in W_v\cap (Q\setminus D)(k)$ for all $v\in T$ such that $C=Q\cap\langle\ell,R\rangle$  is a smooth conic. We  take $P=A$. 

 If $\ell\subset Q$ and $d\ge3$,
we choose a general $R\in(Q\setminus D)(k)$ lying in every
$W_v$ for $v\in T$, and let $C$ be the residual line in
$Q\cap\langle\ell,R\rangle$. Put
$P=\ell\cap C$.
Every $H_i$ is geometrically integral, and $P$ lies on
all of them. Since $\sQ$ is proper, $P$ is integral and
has infinite multiplicity with every boundary component
at every finite place. Thus
$P\in\CC_v(\sQ,\Delta_m)$ for all $v\notin T$.

In either of these two cases, put $F=(D|_C)_{\mathrm{red}}$.
The tuple equal to $R$ at $T$ and to $P$ elsewhere belongs
to $(C,(1-1/m)F)(\A_k)$ and maps into $\mathcal N$.
The proof then follows from Lemma~\ref{lem:conic-reduction}.

%

In the remaining case, where $d=2$ and $\ell\subset Q$,
write $H_i=\ell+C_{s_i}$ for $i=1,2$.
Let $f:Q\to\PP^1$ be the ruling whose fibres are the residual
lines, and put $Z=\{s_1,s_2\}_{\mathrm{red}}$.
We enlarge $T$ to include the places of bad reduction of $Z$.
At each newly added
place $v$, set $W_v=\CC_v(\sQ,\Delta_m)$. The set $f(W_v)$ is open since
$f$ is smooth.
Proposition~\ref{prop:curve} gives a strict $T$-Campana point $t$
for $(\PP^1,(1-1/m)Z)$ such that $t\in f(W_v)$ for every $v\in T$.
By weak approximation on $C_t$, we may find a   rational point $R$
in every $W_v$. The point $P=\ell\cap C_t$ has infinite
multiplicity with $\ell$; outside $T$, its multiplicities
with the other boundary components are those of $t$ with $Z$.
Thus $P$ satisfies the required local conditions for $v\notin T$. Let 
$F=(D|_{C_t})_{\mathrm{red}}=\{P\}$.
The tuple equal to $R$ at $T$ and to $P$ elsewhere belongs
to $(C_t,(1-1/m)F)(\A_k)$ and maps into $\mathcal N$.
The proof then follows from Lemma~\ref{lem:conic-reduction}.
\end{proof}

For the rest of the paper, we may assume that
$E(k)=\varnothing$. Then $q|_W$ is nonzero and nondegenerate, since
the zero form would give $E=L$ with rational points, and
a nonzero radical would give the rational vertex
$\PP(\operatorname{rad}(q|_W))\subset E$.
Thus $E$ is smooth. If it is finite, it is a nonsplit
degree-two point. If it has positive dimension, it is a
smooth geometrically integral quadric.

\section{Quadric surfaces}\label{sec:surface}

Let $Q\subset\PP^3_k$ be a smooth quadric surface, and let
$\Pi_1,\Pi_2$ be arbitrary hyperplanes defined over $k$.
Put
\[
 H_i=Q\cap\Pi_i,\qquad E=H_1\cap H_2,\qquad
 D=H_1+H_2,\qquad \Delta_m=(1-1/m)D_{\mathrm{red}}.
\]
The sections may coincide or share a line. A singular section
is the union of two distinct geometric lines.
Each distinct prime component has weight $m$, with the local
Campana conditions imposed componentwise.

\begin{prop}\label{prop:surface}
For any two hyperplane sections $H_1,H_2$ and every integer
$m\ge2$, $(Q,\Delta_m)$ satisfies Campana strong approximation.
\end{prop}

It suffices to prove the proposition when $E$ is finite.
Indeed, if $\Pi_1\ne\Pi_2$ and $E$ has positive dimension,
then $E=\Pi_1\cap\Pi_2$ is a $k$-line contained in $Q$,
so Proposition~\ref{prop:rational-base-locus} applies.
If $\Pi_1=\Pi_2$, choose an auxiliary hyperplane section $F$
with $H_1\cap F$ finite. Once the finite-intersection case
is proved, applying it to $H_1+F$ and then using monotonicity
\cite[Lemma~2.5(1)]{WeiMarkoff} gives the result for
$D_{\mathrm{red}}=H_1$.

The finite case with $E(k)\ne\varnothing$ also follows from
Proposition~\ref{prop:rational-base-locus}.
For the remaining argument, assume that $E$ is a nonsplit
closed point of degree two. Then $D$ is reduced, so
$\Delta_m=(1-1/m)D$.

\subsection{The conic pencil}

Let $\Lambda=\Pi_1\cap\Pi_2$. The planes through $\Lambda$
give a pencil of conics on $Q$ with base locus $E$.
Let $\pi:X=\operatorname{Bl}_E(Q)\to Q$. The pencil resolves
to a morphism $f:X\to\PP^1_k$. Its boundary
has two vertical fibres, the strict transforms of $H_1$ and $H_2$,
and a horizontal exceptional divisor of degree two over the base.
Every smooth member $C$ other than \(H_1\) and \(H_2\) satisfies
\begin{equation}\label{eq:curve-restrictions}
 H_1|_C=H_2|_C=E.
\end{equation}
More precisely, write $E_{\kbar}=\{A,B\}$. If $H_i$ is singular, its two geometric
lines $L_{i,A}$ and $L_{i,B}$ satisfy
\[
 L_{i,A}|_C=A,\qquad L_{i,B}|_C=B.
\]
A smooth section restricts to $A+B$. Thus the inclusion
$(C,(1-1/m)E)\hookrightarrow(Q,\Delta_m)$ is an orbifold
morphism, and Proposition~\ref{prop:orbifold-continuity}
gives a continuous map on Campana spaces.

Write $K=k(E)$, a quadratic extension of $k$. The nondegenerate
binary quadratic form on $\Lambda$ is a scalar multiple of the
norm form of $K/k$.
We have the normal form
\[
 Q:\quad \N{K}{k}(u)=p(t,w)=at^2+btw+cw^2,
 \qquad H_1=(t=0),\quad H_2=(w=0),
\]
where $u\in K$ represents a pair of $k$-coordinates. Smoothness
of $Q$ gives
\[
 b^2-4ac\ne0.
\]
The sections $H_1$ and $H_2$ are singular exactly when
$c=0$ and $a=0$, respectively.
The member of parameter $\tau=t/w$ is the conic
\begin{equation}\label{eq:fiber}
 C_\tau:\quad \N{K}{k}(u)=p(\tau,1)w^2,
 \qquad t=\tau w.
\end{equation}
Let $Z\subset\PP^1$ be the reduced degree-two scheme
defined by $p(t,w)=0$. It indexes the two singular
geometric fibres. Over a closed point $s\in Z$, the
singular fibre splits over $K.k(s)$.
Define the reduced marked set
\[
 M=\{0,\infty\}\cup Z,\qquad r=\deg M\in\{2,3,4\}.
\]
Here $r=4,3,2$ according as neither, exactly one, or both
boundary sections are singular. 

\subsection{Twisted norm auxiliary varieties}

Assume that $E$ is nonsplit, with splitting field $K$, and use
the reduced marked set $M$ defined above. Write
\[
 Z\setminus\{0,\infty\}=\coprod_j\Spec F_j.
\]
For each $F_j$, choose a quadratic field $A_j/F_j$ containing
a copy of $K$: take $F_jK$ when $K$ does not embed
in $F_j$, and any quadratic extension otherwise. Assign to
$s\in M$ the extensions $L_s/k(s)$ given by $K/k$ at the
boundary parameters and by $A_j/F_j$ at the remaining singular
parameters. If $r\le3$, all marked points are rational and
$L_s=K$ for every $s\in M$.
Each extension $L_s/k(s)$ splits a multiplicity-one component of
$f^{-1}(s)$, since it contains $K$, as permitted in
\cite[Theorem~4.1]{HWW}.

Choose an auxiliary affine coordinate in which $M$ stays away
from infinity, and write $a_s\in k(s)$ for the coordinate of $s$.
For arbitrary $b_s\in k(s)^\times$, let
$p_b:W_b\to\PP^1$ be the standard auxiliary variety of
\cite[Section~3.1]{HWW}, with equations
\[
 \N{L_s}{k(s)}(z_s)=b_s(\lambda-a_s\mu),
 \qquad s\in M,
 \qquad p_b=[\lambda:\mu].
\]
Delete $\lambda=\mu=0$ and the norm-singular loci specified
in that construction. The base map is understood in the
original pencil coordinates after undoing the auxiliary change
of coordinates.

\begin{lem}\label{lem:auxiliary-sa}
For every choice of the nonzero coefficients $b_s$, the variety
$W_b$ satisfies strong approximation off any finite place of $k$.
\end{lem}

\begin{proof}
Over $\kbar$, the norm equations have $r$ blocks:
\[
 z_{j,+}z_{j,-}=\gamma_j(\lambda-a_j\mu),
 \qquad 1\le j\le r,
\]
with $\gamma_j\ne0$ and distinct $a_j$.

Suppose first that $r=4$. Eliminating $\lambda,\mu$ gives
an affine cone over an intersection
$V\subset\PP^7$ of two quadrics. Its pencil has generic rank
eight and exactly four rank-six members. The quadrics have no
common factor. By \cite[Lemmas~1.11 and~1.12]{CTSSD}, $V$ is
geometrically integral and is not a cone.
The Jacobian criterion gives eight singular points.
The standard deleted norm-singular loci have codimension
at least two, and $W_b$ is the remaining open subset of the
smooth cone. Strong approximation then follows from
Theorem~\ref{thm:cone}.

If $r=3$, relabel the singular boundary section as $H_1$.
Then $p(t,w)=a\,t(t-\beta w)$ with $a\beta\ne0$, and the
marked parameters are $0,\infty,\beta$. In these coordinates
the auxiliary equations are
\[
 \N{K}{k}(x)=b_0\lambda,\qquad
 \N{K}{k}(y)=b_\infty\mu,\qquad
 \N{K}{k}(z)=b_\beta(\lambda-\beta\mu).
\]
Eliminating $\lambda,\mu$ gives the nondegenerate rank-six
quadric cone
\[
 \N{K}{k}(z)
 =\frac{b_\beta}{b_0}\N{K}{k}(x)
  -\frac{\beta b_\beta}{b_\infty}\N{K}{k}(y).
\]
Its projectivization is a smooth quadric of dimension four.
The standard deleted loci have codimension at least two,
so \cite[Lemma~6.4]{HWW} gives strong approximation.

If $r=2$, then $p(t,w)=b\,tw$ with $b\ne0$. The equations
are $\N{K}{k}(x)=b_0\lambda$ and
$\N{K}{k}(y)=b_\infty\mu$. Eliminating $\lambda,\mu$
identifies $W_b$ with an open subset of $\Aff^4_k$. 
The standard deleted loci have codimension
at least two, and strong approximation follows from
\cite[Lemma 2.1]{Wei}.
\end{proof}

\subsection{Choosing a fibre and treating new places}


\begin{lem}[Choice of a fibre]\label{lem:choose-conic}
Assume that $E$ is nonsplit. For every nonempty open subset
$\mathcal N\subset(Q,\Delta_m)(\A_k)$, there is
$\tau\in(\PP^1\setminus M)(k)$ such that, for
$j_\tau:C_\tau\hookrightarrow Q$ and
$\delta_\tau=(1-1/m)E$, the set
\[
 j_{\tau,\A}^{-1}(\mathcal N)
       \cap(C_\tau,\delta_\tau)(\A_k)
\]
is nonempty.  Here $j_{\tau,\A}$ denotes the induced adelic map.
\end{lem}

\begin{proof}
Choose $(x_v)\in\mathcal N$. The Hasse principle gives
$Q(k)\ne\varnothing$, so $Q$ is $k$-rational and
$\operatorname{Br}Q=\operatorname{Br}k$.
Apply the descent argument in \cite[Lemma~4.2]{HWW} to
$f:X=\operatorname{Bl}_E(Q)\to\PP^1_k$, with $A=0$ and
$K_s=L_s$, using Lemma~\ref{lem:campana-formal} for its
formal-lemma step. Removing $M$ and the auxiliary
infinity point from the base, together with the exceptional
divisor on $X$, gives an open subset $Q^\circ\subset Q$.
The finite group of cup products used in the descent argument meets
$\operatorname{Br}Q$ only in constants, so orthogonality is
automatic. We obtain $(x_v')\in\mathcal N\cap Q^\circ(\A_k)$,
coefficients $b_s$, and lifts $(z_v)\in W_b(\A_k)$ with
$p_b(z_v)=f(x_v')$. 

Shrink $\mathcal N$ about $(x_v')$ to a basic open set.
Choose a finite set of places $T$ containing the places at which
$\mathcal N$ imposes local restrictions, 
the places ramified in $K/k$, the places of bad reduction of the norm model,
the nonunit places of the defining and twisting coefficients,
and the places where the lifts are not integral. Enlarge $T$ so that the closure of the marked set together
with the auxiliary infinity point is finite \'etale over $\OO_T$.
Outside $T$, use the standard integral model of $W_b$
from \cite[Section~3.5]{HWW}. 
At each place in $T$, choose neighborhoods of $x_v'$ small enough
to retain membership in $\mathcal N$.
Choose a finite place $v_0\notin T$ splitting completely
in the normal closures of all $L_s$.
Lemma~\ref{lem:auxiliary-sa} supplies a rational point on
$W_b$, arbitrarily close to $z_v$ at $T$ and integral
outside $T\cup\{v_0\}$, whose image $\tau$ in $\PP^1_k$ avoids $M$.
Smoothness gives points on $C_\tau$ in the chosen
neighborhoods at $T$.

If $v\notin T\cup\{v_0\}$ and $\tau$ reduces to a marked
point $s$,
the induced place of $k(s)$ has degree one over
$v$, since the marked points have \'etale reduction and
$\tau$ has rational reduction.
By \cite[Proposition~3.8]{HWW}, it extends to a place of
$L_s$ of degree one over $v$.
Since $K$ embeds in $L_s$, the place $v$ splits in $K$.
Splitting at $v_0$ is automatic.
At a split place outside $T$, both points of $E$ are $k_v$-rational.
We choose a local point sufficiently close to either point of $E(k_v)$.
This point can be chosen to have arbitrarily large intersection
multiplicities with the boundary components passing through the
chosen point of $E(k_v)$, and zero intersection multiplicity with the
remaining components. It therefore satisfies the Campana conditions.
At inert places outside $T$, the reduction of
$\tau$ avoids $M$, so $\tau$ and $p(\tau,1)$ are units.
By surjectivity of the unramified norm on units, we may choose
an integral solution $u$ of \eqref{eq:fiber} with $w=1$.
In the point $[u:\tau:1]$, both boundary coordinates are units,
so both boundary multiplicities are zero.
These choices map into $\mathcal N$ and are Campana on
$(C_\tau,\delta_\tau)$ at almost every place, by
\eqref{eq:curve-restrictions}. They therefore give a point
of the required inverse image.
\end{proof}

\begin{proof}[Proof of Proposition~\ref{prop:surface}]
By the reductions above, we may assume that $E$ is a nonsplit
closed point of degree two.
Let $\mathcal N\subset(Q,\Delta_m)(\A_k)$ be nonempty and open. Lemma~\ref{lem:choose-conic} supplies
a smooth conic whose Campana adelic inverse image of
$\mathcal N$ is nonempty.
Lemma~\ref{lem:conic-reduction} gives a strict rational
point in $\mathcal N$. This proves density.
\end{proof}

\begin{cor}\label{cor:surface-one-boundary}
Let $\Sigma$ be a smooth quadric surface and $B$ a hyperplane
section. Then $(\Sigma,(1-1/m)B)$ satisfies Campana strong
approximation for every $m\ge2$.
\end{cor}

\begin{proof}
Apply Proposition~\ref{prop:surface} with $H_1=H_2=B$.
\end{proof}


\section{Reduction to surfaces}\label{sec:higher}

The case $d=2$ was proved in Proposition~\ref{prop:surface}, and the case $E(k)\neq \varnothing$ follows from
Proposition~\ref{prop:rational-base-locus}.
Throughout this section, we assume $d\ge3$ and
$E(k)=\varnothing$.
We reduce the finite and positive-dimensional cases to smooth
quadric surfaces with two and one hyperplane boundary sections,
respectively.

\subsection{Finite common intersection}\label{subsec:finite-intersection}

First assume that $E$ is finite and hence reduced.
Let
\[
 \Lambda_0=\Pi_2\cap\cdots\cap\Pi_d\simeq\PP^2_k,
 \qquad B=Q\cap\Lambda_0,
 \qquad E=B\cap\Pi_1.
\]
The restriction of the defining quadratic form to the line
$\Lambda=\Lambda_0\cap\Pi_1$ has rank two. Hence $B$ is
a reduced plane conic, either smooth or a pair of geometric
lines. In the latter case, its vertex $O$ is $k$-rational and
lies outside $H_1$, since $E$ is reduced.

For the linear algebra below, we write $Q\subset\PP(V)$,
where $V$ is a $k$-vector space of dimension $d+2$ and $q$
is its nondegenerate quadratic form. Choose linear forms
$\ell_i$ defining $\Pi_i$, and put
\[
 V_i=\ker(\ell_i),\qquad
 V_0=\bigcap_{i=2}^d V_i,\qquad W_0=V_0\cap V_1.
\]
Then $\Lambda_0=\PP(V_0)$, $\Pi_1=\PP(V_1)$ and
$\Lambda=\PP(W_0)$, with $\dim V_0=3$ and $\dim W_0=2$.
The reducedness of $E$ makes $q|_{W_0}$ nondegenerate.

\begin{lem}\label{lem:surface-section}
\begin{enumerate}
\item There is a nonempty Zariski-open subset
$\mathcal U\subset Q\setminus D$ such that, for every
$R\in\mathcal U(k)$,
\[
 \Pi_R=\langle\Lambda_0,R\rangle\simeq\PP^3_k,
 \qquad \Sigma_R=Q\cap\Pi_R
\]
is a smooth quadric surface.
\item For every $R$ chosen in \textup{(1)}, the first
boundary section $F_1=H_1|_{\Sigma_R}$ is a reduced plane
conic, and the boundary restrictions satisfy
\[
 H_i|_{\Sigma_R}=B\quad(2\le i\le d),\qquad
 F_1\cap B=E
\]
scheme-theoretically. In particular, $F_1$ and $B$ have
no common component.
\end{enumerate}
\end{lem}

\begin{proof}
\textup{(1)} The parameter space
$\mathcal P=\PP(V/V_0)$ parametrizes four-dimensional
subspaces $U\supset V_0$. The form $q|_{V_0}$ has rank
three or two. In the first case, choose
$z\in V_0^\perp$ with $q(z)\ne0$; in the second, let $o$
span its radical and choose $z$ with $b_q(o,z)\ne0$,
where $b_q$ is the polar form. In either case
$q|_{V_0+kz}$ has rank four. Thus smooth surface sections
form a nonempty open subset of $\mathcal P$.

Projection from $\Lambda_0$ gives
\[
 Q\setminus\Lambda_0\longrightarrow\mathcal P,
 \qquad [r]\longmapsto V_0+kr.
\]
It is dominant, since a general smooth surface section
contains geometric points outside $B$. The inverse image of
the smooth-section locus under this projection, intersected
with $Q\setminus D$,
is the required nonempty open subset $\mathcal U$.

\textup{(2)} For $R=[r]\in\mathcal U(k)$, put
$U_R=V_0+kr$. For $i\ge2$, $\ell_i$ vanishes on $V_0$
and $\ell_i(r)\ne0$, because $R\notin D$.
Hence $V_i\cap U_R=V_0$, giving
$H_i|_{\Sigma_R}=Q\cap\PP(V_0)=B$ and
$F_1\cap B=B\cap\Pi_1=E$ scheme-theoretically.
A hyperplane restriction of the rank-four form defining
the smooth surface $\Sigma_R$ has rank at least two,
so $F_1$ is reduced. Finally, finiteness of $E$ excludes
a common component of $F_1$ and $B$.
\end{proof}

\begin{rmk}\label{rmk:surface-extra-geometry}
	
The section $F_1$ may be singular. Its smoothness defines
an open condition on $R\in Q\setminus D$. This locus is
nonempty by the same nondegeneracy and projection arguments
as in the proof of Lemma~\ref{lem:surface-section}, applied
to plane sections of $Q\cap\Pi_1$ through $\Lambda$. Indeed,
$q|_{W_0}$ is nondegenerate and
$\operatorname{rank}(q|_{V_1})\ge d\ge3$.
Denote this nonempty open subset by $\mathcal U_1$.
Since $Q$ is irreducible, $\mathcal U\cap\mathcal U_1$
is nonempty and open. Thus the smoothness of $F_1$ may
also be imposed if desired.	
	The reducedness of $E=F_1\cap B$ forces the vertex of
either singular conic to lie outside the other. Both conics are
therefore smooth at $E$ and meet transversely there. 
The extra smoothness of $F_1$ and the separate transversality
assertion are not required for Proposition~\ref{prop:surface} to apply.
\end{rmk}

The following lemma gives the nonempty adelic set needed
on the chosen surface.

\begin{lem}\label{lem:new-surface-places}
Let $\mathcal N\subset(Q,\Delta_m)(\A_k)$ be nonempty and
open. There exists $R\in\mathcal U(k)$ such that, with
$\Sigma=\Sigma_R$, $\delta_R=(1-1/m)(F_1+B)$ and
$j:\Sigma\hookrightarrow Q$, the set
\[
 \mathcal G=j^{-1}_{\A}(\mathcal N)
       \cap(\Sigma,\delta_R)(\A_k)
\]
is nonempty.
\end{lem}

\begin{proof}
Shrink $\mathcal N$ to a basic open subset, taking its finite
set $T$ of local conditions large enough to contain the model
and reduction exceptions for $B$ and $E$.
The Hasse principle and weak approximation give
$R\in\mathcal U(k)$ satisfying the local conditions at $T$.

If $v\notin T$, we lift a residue-field point outside $E$ when $B$ is
smooth, and otherwise take its vertex. Hence we may choose $b_v\in B(k_v)$ reducing outside
$H_1$. Since
these points lie on $H_i$ for $i\ge2$ and have zero multiplicity
with $H_1$, they satisfy the required Campana conditions on $Q$.
The tuple equal to $R$ at $T$ and to $b_v$ elsewhere belongs
to $(\Sigma,\delta_R)(\A_k)$, since it lies
on $B$ and has zero multiplicity with $F_1$ at almost every place.
Its image lies in $\mathcal N$, proving the assertion.
\end{proof}

\begin{rmk}\label{rmk:weaker-hypotheses}
Finiteness of $E$ implies that the hyperplanes are linearly
independent and that $\deg E=2$.
Our assumption $E(k)=\varnothing$ makes $E$ reduced, which
forces the intermediate plane conic $B$ to be reduced,
although it may be singular.
A nonreduced $E$ would have the form $2P$ with $P\in Q(k)$,
a case already covered by
Proposition~\ref{prop:rational-base-locus}. In this case, $B$ may even
be a double line.
Thus no separate reducedness or smoothness hypothesis is
needed in the finite-intersection argument.
\end{rmk}

\subsection{Positive-dimensional common intersection}\label{subsec:positive-intersection}

Now assume that $\dim E\ge1$.
With $L=\PP(W)$ as in Section~\ref{sec:rational-base-locus},
the form $q|_W$ is anisotropic over $k$. Choose a $k$-plane
$L'=\PP(V')\subset L$ and put $B=Q\cap L'$. The restriction
$q|_{V'}$ is nondegenerate, since a nonzero radical vector
would be isotropic. Thus $B$ is a smooth conic.

The rank-three argument in the proof of
Lemma~\ref{lem:surface-section}\textup{(1)} gives a nonempty
open subset $\mathcal U_+\subset Q\setminus D$ such that
\[
 \Sigma_R=Q\cap\langle L',R\rangle
\]
is a smooth quadric surface for every $R\in\mathcal U_+(k)$.
Since $R\notin D$ and $L'\subset\Pi_i$, we have
\[
 H_i|_{\Sigma_R}=B\qquad(1\le i\le d).
\]
Thus the reduced boundary on $\Sigma_R$ is the smooth conic $B$.

\begin{proof}[Proof of Theorem~\ref{thm:main}]
Let $\mathcal N\subset(Q,\Delta_m)(\A_k)$ be nonempty and open.

For finite $E$, Lemma~\ref{lem:new-surface-places} gives a
smooth surface $\Sigma$ with boundary
$\delta=(1-1/m)(F_1+B)$ and a nonempty set
$\mathcal G=j_{\mathbf A}^{-1}(\mathcal N)\cap(\Sigma,\delta )(\A_k)$.
For positive-dimensional $E$, we choose a basic subset of
$\mathcal N$ with a finite set $T$ of local conditions,
large enough that $B(k_v)\ne\varnothing$ outside $T$.
By the Hasse principle and weak approximation, we choose
$R\in\mathcal U_+(k)$ satisfying these local conditions.
Put $\Sigma=\Sigma_R$ and $\delta=(1-1/m)B$.
The tuple equal to $R$ at $T$ and to points of $B(k_v)$
elsewhere belongs to $(\Sigma,\delta)(\A_k)$ and maps into
$\mathcal N$.
Thus $\mathcal G$ is again nonempty.

In both cases $j:(\Sigma,\delta)\hookrightarrow(Q,\Delta_m)$
is an orbifold morphism, so $\mathcal G$ is open by
Proposition~\ref{prop:orbifold-continuity}.
Applying Proposition~\ref{prop:surface} in the finite case and
Corollary~\ref{cor:surface-one-boundary} in the positive-dimensional
case gives a strict rational point in $\mathcal G$.
Its image is strict on $Q$ and belongs to $\mathcal N$,
proving density.
\end{proof}

\begin{rmk}\label{rm:geometric-boundary}
More generally, suppose $D$ is an effective $k$-divisor such that
		\[
		D_{\kbar}=\sum_{i=1}^d H_i,\qquad
		H_i=Q_{\kbar}\cap\Pi_i.
		\]
		We assume that the sum $\sum_{i=1}^d \Pi_i$ is defined over $k$,
		although each $\Pi_i$ need not be defined over $k$. Let $L=\bigcap_i\Pi_i$ and $E=Q\cap L$;  these are defined over $k$.
		For $m\ge2$, let $\Delta_m=(1-1/m)D_{\mathrm{red}}$.
	Then $(Q,\Delta_m)$ satisfies Campana strong approximation if $d=2$, or if $d\ge3$ and at least one of the following holds: $E(k)\ne\varnothing$, $\dim E>0$, or some $\Pi_i$ is defined over $k$. 
		
		If $E(k)\ne\varnothing$, or if $d\geq 3$ and either some $\Pi_i$ is defined over $k$ or $\dim E>0$,  the proof of
		Proposition~\ref{prop:rational-base-locus} and  reductions
		in Section~\ref{subsec:finite-intersection} and \ref{subsec:positive-intersection}  apply
		with the same boundary restrictions.
		On surfaces, it remains to consider that both hyperplanes are conjugate. Their boundary parameters then form a closed point of degree two. The marked set has degree four when the sections are smooth and degree two when they are singular. The corresponding auxiliary-variety arguments therefore apply with minor modifications.
\end{rmk}

\end{document}